\documentclass[12pt]{article}  

\usepackage{graphicx}

\usepackage{amsmath, amssymb, amsfonts}
\usepackage{color}
\usepackage{longtable}
\usepackage[justification=centering]{caption}

\usepackage{tikz,pgf}

\usepackage[colorlinks=true,urlcolor=blue,bookmarksopen=false,pdfpagemode=UseNone,pdfstartview=FitH]{hyperref}
\usepackage[all]{hypcap}

\newtheorem{Thm}{Theorem}[section]
\newtheorem{Lem}[Thm]{Lemma}
\newtheorem{Prop}[Thm]{Proposition}

\newtheorem{Ques}[Thm]{Question}
\newtheorem{BFST}{\bfseries Theorem~\ref{mucpp} \mdseries (\cite{BFS11})}

\newcommand{\qed}{\rule{5pt}{5pt}}
\newcommand{\prf}{{\it Proof.\hspace{5pt}}}
\newcommand{\bbP}{{\mathbb P}}
\DeclareMathOperator{\st}{st}
\newcommand{\bea}{\begin{eqnarray*}}
\newcommand{\eea}{\end{eqnarray*}}

\newcommand{\bfig}[1]{\begin{figure}[#1] \begin{center} \addtocounter{Thm}{1}
}
\newcommand{\efig}{\end{center} \end{figure}\vskip0\baselineskip}

\newcommand{\btab}[1]{\begin{table}[#1] \begin{center} \addtocounter{Thm}{1}
}
\newcommand{\etab}{\end{center} \end{table}\vskip0\baselineskip}

\definecolor{beamerblue}{rgb}{.2,.2,.7}

\definecolor{brightbblue}{rgb}{.6,.7,1}

\title{Discrete Morse theory and the consecutive pattern poset
}

\author{
Bruce E. Sagan\\[-5pt]
\small Department of Mathematics, Michigan State University,\\[-5pt]
\small East Lansing, MI 48824-1027, USA\\[-5pt]
\small \url{http://www.math.msu.edu/~sagan}\\[-5pt]
\small \href{mailto:sagan@math.msu.edu}{\texttt{sagan@math.msu.edu}}\\[6pt]
Robert Willenbring\\[-5pt]
\small Department of Mathematics, University of Mary,\\[-5pt]
\small Bismarck, North Dakota, 58504\\[-5pt]
\small \href{mailto:willenb4@msu.edu}{\texttt{willenb4@msu.edu}}
}

\date{\today\\[10pt]
	\begin{flushleft}
	\end{flushleft}
}

\begin{document}
\maketitle

\begin{abstract} We use discrete Morse theory to provide another proof of Bernini, Ferrari, and Steingr\'imson's formula for the M\"obius function of the consecutive pattern poset.  In addition, we are able to determine the homotopy type of this poset.  Earlier, Bj\"orner determined the M\"obius function and homotopy type of factor order and the results are remarkably similar to those in the pattern case.  In his thesis, Willenbring used discrete Morse theory to give an illuminating proof of Bj\"orner's result.  Since our proof parallels Willenbring's, we also consider the relationship between the two posets.  In particular, we show that some of their intervals are isomorphic, and also that there is a sequence of posets interpolating between the two all of whom have essentially the same M\"obius function.
\end{abstract}
\maketitle

\section{Introduction}
\label{intro}

The M\"obius function of the consecutive pattern poset was determined by Bernini, Ferrari, and Steingr\'imson~\cite{BFS11}.  We give another proof of their result using discrete Morse theory.   In our demonstration, the definitions needed to state this formula naturally arise from the structure of the poset.  Furthermore, the use of discrete Morse theory allows us to easily determine the homotopy type of the poset.  Bj\"orner's formula for the M\"obius function of factor order~\cite{aB93} is remarkably similar to that of the consecutive pattern poset, so we will also consider the relationship between these two posets. 

Let $S_d$ be the set of all permutations of the first $d$ positive integers. A \emph{consecutive} pattern $\sigma=\sigma(1)\sigma(2)\ldots \sigma(k)\in S_k$ appears in a permutation $\tau=\tau(1)\tau(2)\ldots \tau(d) \in S_d$ if the letters of some subsequence $\tau(i+1)\tau(i+2)\ldots \tau(i+k)$ of $\tau$ appear in the same order of size as the letters in $\sigma$.  The consecutive pattern poset is $S=\cup_{d>0}S_d$ ordered with respect to consecutive pattern containment.  The \emph{length} of $\tau$, denoted $|\tau|$, is the number of positions in $\tau$.  In particular, if $\tau\in S_d$, $|\tau|=d$.  A permutation $\sigma$ is \emph{monotone} if its letters are either strictly increasing or strictly decreasing.  The \emph{standard form} of a sequence of distinct integers $s=s(1)s(2)\ldots s(d)$ is the permutation $p=p(1)p(2)\ldots p(d)$ in $S_d$ whose letters appear in the same order of size as the letters in $s$. 

To state the formula for the M\"obius function of the consecutive pattern poset as a simple recursive formula, we make the following definitions.  A \emph{prefix of length $k$} of a permutation $\tau\in S_d$  is the standard form of the first $k$ letters of $\tau$.  If $k<d$, then the prefix is \emph{proper}.  For example, prefixes of $53412$ include $312$ and $4231$.  Similarly, a \emph{suffix of length $k$} of a permutation $\tau\in S_d$ is the standard form of the last $k$ letters of $\tau$.  We define the \emph{interior} of a permutation $\tau=\tau(1)\tau(2)\ldots \tau(d)\in S_d$ for $d>2$ to be the standard form of the sequence $\tau(2)\ldots \tau(d-1)$.  The interior of $\tau$ will be denoted $i(\tau)$.  We define the \emph{exterior} of a permutation $\tau$ to be the longest permutation which is the standard form of both a proper prefix and a suffix $\tau$. The exterior of $\tau$ will be denoted $x(\tau)$.  Notice the exterior is always defined, since the standardization of a single integer is 1.

The following theorem is an equivalent statement of a theorem of Bernini, Ferrari, and Steingr\'imson.  For those familiar with this work, we note that we use the condition $x(\tau)\not\leq i(\tau)$ in place of the carrier element of $\tau$. 

\begin{Thm}[\cite{BFS11}]
\label{mucpp}
In the consecutive pattern poset $S$, if $\sigma\leq \tau$ then 
$$\mu(\sigma,\tau)=\begin{cases}
\mu(\sigma,x(\tau)) & \mbox{if } |\tau|-|\sigma|>2\mbox{ and } \sigma\leq x(\tau)\not\leq i(\tau),\\
1 &\mbox{if } |\tau|-|\sigma|=2\mbox{, } \tau \mbox{ is not monotone, and }\sigma\in\{i(\tau),x(\tau)\},\\
(-1)^{|\tau|-|\sigma|} &\mbox{if } |\tau|-|\sigma|<2,\\
0 &\mbox{otherwise.}\end{cases}$$
\end{Thm}

To state our proof of this result, we will need to use a theorem of Babson and Hersh~\cite{BH05}. This result gives a way of applying Forman's discrete version of Morse theory~\cite{rF95} to partially ordered sets.  For brevity, we are only going to state the minimum number of definitions to apply this theorem.  A reader who is interested in learning more about discrete Morse theory is encouraged to begin with Robin Forman's outstanding introduction the topic~\cite{rF02}.

Let $P$ be any poset.  The notation $x\rightarrow y$ will indicate that $x$ covers $y$ in $P$.  Let $C: x=z_0 \stackrel{ }{\rightarrow}z_1\stackrel{ }{\rightarrow}\ldots\stackrel{ }{\rightarrow}z_n=y$ be a chain in $P$.  Since each pair of adjacent elements are related by a cover, $C$ is called a \emph{saturated chain}.  The \emph{closed interval} of $C$ from $z_i$ to $z_j$ is the chain $C[z_i,z_j]: z_i \stackrel{ }{\rightarrow}z_{i+1}\stackrel{ }{\rightarrow}\ldots\stackrel{ }{\rightarrow}z_j$.  The \emph{open interval} of $C$ from $z_i$ to $z_j$, $C(z_i,z_j)$, is defined similarly.  The closed interval $C[z_{i},z_{i}]$ consisting of the single element $z_i$ will also be written $z_i$, but the context will always indicate whether we are referring to the element or the interval.  Notice that the interval $[x,y]$ is non-empty when $x\leq y$ in the poset $P$, while $C[z_i,z_j]$ is non-empty when $z_i\geq z_j$.  A  saturated chain $C$ of the interval $[x,y]$ is a \emph{maximal chain} if $z_0=y$ and $z_n=x$.

Given two maximal chains  $C: z_0 \stackrel{ }{\rightarrow}z_1\stackrel{ }{\rightarrow}\ldots\stackrel{ }{\rightarrow}z_n$ and  $D: x_0 \stackrel{ }{\rightarrow}x_1\stackrel{ }{\rightarrow}\ldots\stackrel{ }{\rightarrow}x_n$ in an interval $[x,y]$, we say $C$ and $D$ \emph{agree to index $k$} if $z_i=x_i$ for all $i\leq k$.  We say $C$ and $D$ \emph{diverge from index $k$} if $C$ and $D$ agree to index $k$ and $z_{k+1}\neq x_{k+1}$.  A total ordering $C_1<C_2<\ldots<C_n$ of the maximal chains of an interval is a \emph{poset lexicographic order} if it satisfies the following: suppose $C<D$ and $C$ and $D$ diverge from index $k$; if $C'$ and $D'$ agree to index $k+1$ with $C$ and $D$, respectively, then $C'<D'$.

Suppose $C_1<C_2<\ldots<C_n$ is an ordering of the maximal chains of the closed interval $[x,y]$. An interval $C(z_i,z_j)$ is a \emph{skipped interval} of a maximal chain $C$ if
$$C-C(z_i,z_j) \subseteq C' \mbox{ for some } C' < C.$$
It is a \emph{minimal skipped interval (MSI)} if it does not properly contain another skipped interval.  We write $I(C)$ for the set of all MSIs of a chain $C$.  To find the set $I(C)$, consider each interval $I\subseteq C(y,x)$ and see if $C-I\subseteq C'$ for any $C'\subseteq C$, then throw out any such interval that is nonminimal.  In Table~\ref{table:JC}, we give examples of minimal skipped intervals in the context of the consecutive pattern poset.

Notice $I(C)$ could contain intervals which overlap, that is, intervals with non-empty intersection. Babson and Hersh's result requires a set of disjoint intervals derived from $I(C)$, which we will denote $J(C)$.  We construct $J(C)=\{J_1,J_2,\ldots\}$ as follows.  Order the intervals of $I(C)$ based on when they are first encountered in $C$.  Thus, $I_1$ will contain the element $z_i$ of smallest index that appears in any interval in $I(C)$, $I_2$ will contain the element $z_j$ of smallest index that appears in any interval in $I(C)$ other than $I_1$, etc.  Let $J_1=I_1$.  Then consider the intervals $I'_2=I_2-J_1$, $I'_3=I_3-J_1$, and so forth. Throw out any that are not containment minimal, and pick the first one that remains to be $J_2$.  Continue this process until no intervals remain to add to $J(C)$.

The set of intervals $J(C)$ \emph{covers} $C$ if its union equals the open interval $C(z_0,z_n)$.  A chain $C$ is called \emph{critical} if $J(C)$ covers $C$.  Finally, if a chain $C$ is critical, the \emph{critical dimension} of the chain is 
$$
d(C)=\#J(C)-1
$$
where $\#$ denotes cardinality.

\begin{Thm}[\cite{BH05}]
\label{chainmu}
For any poset lexicographic order on the maximal chains of $[x,y]$,
$$\mu(x,y)=\sum_{C}(-1)^{d(C)},$$
where the sum is over all critical chains $C$ in the poset lexicographic order.\hfill\qed
\end{Thm} 

The rest of the paper is organized as follows.  In Section~\ref{cpp}, we use Theorem~\ref{chainmu} to reprove Bernini, Ferrari, and Steingr\'imson's formula in Theorem~\ref{mucpp}.  In Section~\ref{cppofo}, we define factor order and consider the relationship between posets ordered by factor order and the consecutive pattern poset.  In particular, we show that certain intervals in $S$ are isomporphic to intervals in factor order using an alphabet of two letters.  We also define a sequence of posets interpolating between $S$ and factor order on the integers each of which has a M\"obius function which can be described as in Theorem~\ref{mucpp}.

\section{Discrete Morse Theory and the Consecutive Pattern Poset}
\label{cpp}

To get a sense of the structure of the consecutive pattern poset, we first consider the covering relations.  The following lemma is easy to prove and so its demonstration is left to the reader.

\begin{Lem}[\cite{BFS11}]
\label{cover}
A permutation $\tau=\tau(1)\tau(2)\ldots \tau(d)$ in the consecutive pattern poset can only cover the standard form of the sequences $\tau(2)\ldots \tau(d)$ and $\tau(1)\ldots \tau(d-1)$.  These two permutations are distinct unless $\tau$ is monotone, in which case they are equal and monotone.\hfill\qed
\end{Lem}

An \emph{expansion} $\eta$ of $\sigma \in S_k$ is $\eta$ of the form $0\ldots 0\sigma 0\ldots 0$.   An \emph{embedding} of $\sigma$ into $\tau=\tau(1)\tau(2)\ldots \tau(d)$ is an expansion $\eta$ of $\sigma$ with length $|\tau|$ such that the nonzero entries of $\eta$, $\eta({i+1}) \ldots \eta({i+k})$, are the standard form of the sequence $\tau(i+1)\tau(i+2)\ldots \tau(i+k)$.  In addition, we say each entry $\tau(1),\ldots, \tau(i)$ and $\tau(i+k+1),\ldots, \tau(d)$ is \emph{reduced to $0$} in $\eta$. Note that when $\tau$ covers $\sigma$, there is a unique embedding unless $\tau$ is monotone, in which case we will always be using the suffix embedding $0\sigma$.

\btab{t!}
\begin{center} $I(C)$ intervals for for $[1,213546]$:\end{center}

\begin{longtable}[!h]{@{\extracolsep{-5pt}}rccccccccccc}
Chain Id & $\rho_0$ & $l_1$ & $\rho_1$  & $l_2$ & $\rho_2$  & $l_3$ & $\rho_3$  & $l_4$ & $\rho_4$  & $l_5$ & $\rho_5$ \\
1-2-3-4-5  & $213546$ & $1$  & $012435$ & $2$  & $001324$ & $3$  & $000213$ & $4$  & $000012$ & $5$  & $000001$ \\
1-2-3-6[-]4 & $213546$ & $1$ & $012435$ & $2$ & $001324$ & $3$ & $000213$ & $6$ & $[000210]$ & $4$ & $000010$ \\
1-2-6[-]3-4 & $213546$ & $1$ & $012435$ & $2$ & $001324$ & $6$ & $[001320]$ & $3$ & $000210$ & $4$ & $000010$ \\
1-2-6[-]5[-]3 & $213546$ & $1$ & $012435$ & $2$ & $001324$ & $6$ & $[001320]$ & $5$ & $[001200]$ & $3$ & $000100$ \\
1-6[-]2-3-4 & $213546$ & $1$ & $012435$ & $6$ & $[012430]$ & $2$ & $001320$ & $3$ & $000210$ & $4$ & $000010$ \\
1-6[-]2-5[-]3 & $213546$ & $1$ & $012435$ & $6$ & $[012430]$ & $2$ & $001320$ & $5$ & $[001200]$ & $3$ & $000100$ \\
1-6-5[-]2-3 & $213546$ & $1$ & $012435$ & $6$ & $012430$ & $5$ & $[012300]$ & $2$ & $001200$ & $3$ & $000100$ \\
6[-]1-2-3-4 & $213546$ & $6$ & $[213540]$ & $1$ & $012430$ & $2$ & $001320$ & $3$ & $000210$ & $4$ & $000010$ \\
6[-]1-2-5[-]3 & $213546$ & $6$ & $[213540]$ & $1$ & $012430$ & $2$ & $001320$ & $5$ & $[001200]$ & $3$ & $000100$ \\
6[-]1-5[-]2-3 & $213546$ & $6$ & $[213540]$ & $1$ & $012430$ & $5$ & $[012300]$ & $2$ & $001200$ & $3$ & $000100$ \\
6-5[-]1-2-3 & $213546$ & $6$ & $213540$ & $5$ & $[213400]$ & $1$ & $012300$ & $2$ & $001200$ & $3$ & $000100$ \\
6[-5-]4[-]1-2 & $213546$ & $6$ & $[213540$ & $5$ & $213400]$ & $4$ & $[213000]$ & $1$ & $012000$ & $2$ & $001000$ \\
6[-5[-]4-]3[-]1 & $213546$ & $6$ & $[213540$ & $5$ & $[213400]$ & $4$ & $213000]$ & $3$ & $[210000]$ & $1$ & $010000$ \\
\end{longtable}

\begin{center} $J(C)$ intervals for $[1,213546]$\end{center}
\begin{center} where $I(C)=J(C)$ except for the last chain: \end{center}

\begin{longtable}[!h]{@{\extracolsep{-5pt}}rccccccccccc}
Chain Id & $\rho_0$ & $l_1$ & $\rho_1$  & $l_2$ & $\rho_2$  & $l_3$ & $\rho_3$  & $l_4$ & $\rho_4$  & $l_5$ & $\rho_5$ \\
6[-5-]4[-]3[-]1 & $213546$ & $6$ & $[213540$ & $5$ & $213400]$ & $4$ & $[213000]$ & $3$ & $[210000]$ & $1$ & $010000$\\
\end{longtable}
\caption{Comparing $I(C)$ and $J(C)$ for our poset lexicographic ordering of the maximal chains of $[1,213546]$.}
\label{table:JC}
\etab

Let $[\sigma,\tau]$ be an interval in $\cup_{d>0}S_d$.  Let 
$$C:\tau=\rho_0\stackrel{l_1}{\rightarrow}\rho_1\stackrel{l_2}{\rightarrow}\ldots\stackrel{l_{n-1}}{\rightarrow} \rho_{n-1}\stackrel{l_n}{\rightarrow}\rho_n=\sigma$$
be a maximal chain in $[\sigma,\tau]$, where the $l_i$ are defined as follows.  Using the convention of the previous paragraph for the monotone case, each cover $\rho_{i-1}\rightarrow\rho_i$ defines a unique embedding of $\rho_i$ in $\rho_{i-1}$.  Inductively, this defines an embedding $\eta_i$ of $\rho_i$ in $\tau$.  We let $l_i$ be the position zeroed out in passing from $\eta_{i-1}$ to $\eta_i$.  Examples of maximal chains written in this form are given in Table~\ref{table:JC}.

Note that this construction gives each maximal chain its own unique sequence $l_1\ldots l_n$ which we can use to identify it.  We call this sequence a maximal chain's \emph{chain id}.  It is easy to see that lexicographically ordering these chain ids produces a poset lexicographic order on the maximal chains of $[\sigma,\tau]$.  This is the order we will use to find the MSIs of the maximal chains $C$, and ultimately the sets $J(C)$ which will allow us to apply Theorem~\ref{chainmu}.

Table~\ref{table:JC} also contains an example of this poset lexicographic ordering and the resulting MSIs.  The table signifies the MSIs of each maximal chain by surrounding the appropriate permutations with brackets, and contains one example of a maximal chain with overlapping MSIs, namely, the last chain.

To facilitate the exposition, we make the following definitions. A \emph{descent} in a maximal chain is a permutation $\rho_i$ where $l_i>l_{i+1}$.  We say $\rho_i$ is a \emph{strong descent} if $l_i>l_{i+1}+1$, and a \emph{weak descent} if $l_i=l_{i+1}+1$.  An \emph{ascent} in a maximal chain is a permutation $\rho_i$ where $l_i<l_{i+1}$.  In Table~\ref{table:JC}, the chain with chain id $6-5-1-2-3$ has $\rho_1=21354$ as a weak descent, $\rho_2=2134$ as a strong descent, and $\rho_3=123$ as an ascent.

\begin{Lem}
\label{descent}
Let $C:\tau=\rho_0\stackrel{l_1}{\rightarrow}\rho_1\stackrel{l_2}{\rightarrow}\ldots\stackrel{l_{n-1}}{\rightarrow} \rho_{n-1}\stackrel{l_n}{\rightarrow}\rho_n=\sigma$ be a maximal chain in $[\sigma,\tau]$.  If $\rho_i$ is a strong descent, then $\rho_i$ is an MSI in $C$.
\end{Lem}

\prf
 Since $l_{i}-l_{i+1}>1$, we know that  $\rho_i$ is the standard form of the longest proper prefix of $\rho_{i-1}$ and $\rho_{i+1}$ is the standard form of the interior of $\rho_{i-1}$.  Let $u$ be the longest proper suffix of $\rho_{i-1}$. 
 
Suppose $u$ is not a monotone permutation.  Then
$$C':\rho_0\stackrel{l_1}{\rightarrow}\ldots\stackrel{l_2}{\rightarrow} \rho_{i-2}\stackrel{l_{i-1}}{\rightarrow} \rho_{i-1}\stackrel{l_{i+1}}{\rightarrow}u\stackrel{l_{i}}{\rightarrow}\rho_{i+1}\stackrel{l_{i+2}}{\rightarrow}\ldots\stackrel{l_{n}}{\rightarrow} \rho_n$$ 
is a lexicographically earlier chain than $C$.  Hence, $\rho_i$ is a skipped interval in $C$.  Since $\rho_i$ is an interval consisting of a single element, $\rho_i$ is an MSI.

Suppose $u$ is a monotone permutation.  Then $l_i$ cannot be reduced in $u$.  However, the chain
$$D':\rho_0\stackrel{l_1}{\rightarrow}\ldots\stackrel{l_2}{\rightarrow} \rho_{i-2}\stackrel{l_{i-1}}{\rightarrow} \rho_{i-1}\stackrel{l_{i+1}}{\rightarrow}u\stackrel{l_{i+2}}{\rightarrow}\rho_{i+1}\stackrel{l_{i+3}}{\rightarrow}\ldots\stackrel{l_{n}}{\rightarrow} \rho_{n-1}\stackrel{l_{n}+1}{\rightarrow} \rho_{n}$$
is a lexicographically earlier chain than $C$ because each $\rho_j$ for $j>i$ is monotone.  Hence, $\rho_i$ is an MSI in this case as well.
\hfill\qed

Considering this result in Table~\ref{table:JC} suggests that Lemma~\ref{descent} accounts for a small proportion of the MSIs in this poset lexicographic order.  The next result, however, gives a great deal of information about the critical chains in this order.

\begin{Lem}
\label{ascent}
Let $C:\tau=\rho_0\stackrel{l_1}{\rightarrow}\rho_1\stackrel{l_2}{\rightarrow}\ldots\stackrel{l_{n-1}}{\rightarrow} \rho_{n-1}\stackrel{l_n}{\rightarrow}\rho_n=\sigma$ be a maximal chain in $[\sigma,\tau]$.  If $\rho_i$ is an ascent, then it is not contained in any MSI.
\end{Lem}

\prf
We will prove this lemma by contradiction.  Suppose $C[\rho_r,\rho_s]$ is an MSI that contains $\rho_i$. Notice that $\rho_i$ may only be preceded by ascents in this interval because if there are descents, the last one that occurs before $\rho_i$ would be a strong descent.  By Lemma~\ref{descent}, this would be an MSI, contradicting the minimality of $C[\rho_r,\rho_s]$.  Thus, it suffices to derive a contradiction for $\rho_i=\rho_{r}$, the first ascent in the interval $C[\rho_r,\rho_s]$.

Since $C[\rho_r,\rho_s]$ is an MSI of $C[\tau,\sigma]$ if and only if it is an MSI of $C[\rho_{r-1},\rho_{s+1}]$, it suffices to consider the case $r=1$.   However, if $r=1$ then $\rho_1$ being an ascent forces $l_1=1$.  This implies $\rho_1$ appears in all chains preceding $C$.  Therefore, $\rho_1$ can be removed from any skipped interval in which it appears and that interval will still be skipped, contradicting the fact that $C[\rho_1,\rho_s]$ is minimal.
\hfill\qed

Notice that Lemma~\ref{ascent} implies that all MSIs of a chain $C$ consist entirely of descents.  Thus, only the lexicographically last chain in an interval can possibly be critical.  The next lemma covers the two basic cases of MSIs in a chain $C$.  We have already encountered the first case, while the second case is new. 

\begin{Lem}
\label{inout} Suppose $\tau$ is not monotone and $|\tau|\geq2$:\begin{enumerate}
\item There are two maximal chains in the interval $[i(\tau),\tau]$, and if $C$ is the second chain, then it has a decreasing chain id and a unique MSI, $C(\tau,i(\tau))$.
\item If $x(\tau)\not\leq i(\tau)$, then there are two maximal chains in the interval $[x(\tau),\tau]$.  If $C$ is the second chain, then it has a decreasing chain id and a unique MSI, $C(\tau,x(\tau))$.
\end{enumerate}
\end{Lem}

\prf
The first case follows from Lemma~\ref{descent}.  For the second case,   once we remove the first or last element, there is only one copy of $x(\tau)$ left in $\tau$.  Therefore, there are two maximal chains in the interval $[x(\tau),\tau]$: the first chain, which ends at the suffix embedding, and the last chain, which ends at the prefix embedding.  Since $x(\tau)$ is the longest prefix of $\tau$ which is also a suffix, and the first chain contains only suffixes of $\tau$ while the last chain contains only prefixes of $\tau$, these two chains have only $x(\tau)$ and $\tau$ in common.  Therefore, $C(\tau,x(\tau))$ is an MSI, completing the proof.
\hfill\qed

We can illustrate this lemma using the intervals $[123,21354]$ and $[21,21354]$.  Note the interior of $21354$ is $123$, so that the maximal chains of $[123,21354]$ are $21354\rightarrow 1243\rightarrow 123$ and $21354\rightarrow 2134\rightarrow 123$, giving $C(21354,123)$ as an MSI in the second chain.  The exterior  element of $21354$ is $21$ and the interior element is $123$ with $21\not\leq 123$, so  the maximal chains of $[21,21354]$ are $21354\rightarrow 1243\rightarrow 132\rightarrow 21$ and $21354\rightarrow 2134\rightarrow 213\rightarrow 21$, giving $C(21354,21)$ as an MSI in the second chain.

Note there can be no overlap between these two types of MSI.  Proposition~\ref{outerword} generalizes Lemma~\ref{inout}(2), and the theorem that follows shows that we have identified all cases of MSIs.  It will be convenient to adopt the convention that a sequence $l_{i+1}$ consisting of a single label is not decreasing, corresponding to the fact that the interval $C(\rho_i,\rho_{i+1})$ is empty and so not an MSI.

\begin{Prop}
\label{outerword}
Let $C:\tau=\rho_0\stackrel{l_1}{\rightarrow}\rho_1\stackrel{l_2}{\rightarrow}\ldots\stackrel{l_{n-1}}{\rightarrow} \rho_{n-1}\stackrel{l_n}{\rightarrow}\rho_n=\sigma$ be a maximal chain in the interval $[\sigma,\tau]$. Suppose there are $i$ and $j$ such that $\rho_j=x(\rho_i)\not\leq i(\rho_i)$, and the sequence $l_{i+1},\ldots, l_j$ is decreasing.  Then $C(\rho_{i},\rho_{j})$ is an MSI in $C$.
\end{Prop}

\prf
Since the sequence $l_{i+1},\ldots, l_j$ is decreasing, $\rho_i$ can not be monotone and $j\geq i+2$.  So Lemma~\ref{inout}(2) implies that $C(\rho_{i},\rho_{j})$ is an MSI in the subchain of $C$ that is its intersection with $[\rho_j,\rho_i]$.  The proof is concluded by noting that if $\rho\leq \phi$ are elements of any maximal chain $C$ then $C(\phi,\rho)$ is an MSI of $C[\phi,\rho]$ if and only if it is an MSI of $C$.
\hfill\qed

\begin{Thm}
\label{MSIchar}
The interval $C(\rho_i,\rho_j)$ is an MSI of a maximal chain $C$ if and only if $C(\rho_i,\rho_j)=\rho_{i+1}$ and $\rho_{i+1}$ is a strong descent, or $\rho_j=x(\rho_i)\not\leq i(\rho_i)$, and the sequence $l_{i+1},\ldots, l_j$ is decreasing.
\end{Thm}

\prf
The reverse implication follows from Lemma~\ref{descent} and Proposition~\ref{outerword}.

Suppose $C(\rho_i,\rho_j)$ is an MSI in $C$.  By Lemma~\ref{ascent}, the sequence $l_{i+1},\ldots,l_{j}$ is a decreasing sequence. If $\rho_{i+1}$ is a strong descent, then $\rho_{i+1}$ is an MSI by Lemma~\ref{descent}.  This implies $C(\rho_i,\rho_j)=\rho_{i+1}$.  If $\rho_{i+1}$ is not a strong descent then, by containment minimality, none of the descents are strong descents.  Also, our sequence is decreasing, so we conclude that $\rho_j$ is a prefix of the permutation $\rho_i$.  

Since there is only one chain ending at the prefix embedding of $\rho_j$ in $\rho_i$, and $[\rho_i,\rho_j]$ contains another chain previous to $C(\rho_j,\rho_i)$, there must be a second copy of $\rho_j$ in $\rho_i$.
 Let $k$ be the largest index so that $\rho_k$ contains at least two copies of $\rho_j$.  Then $\rho_{k+1}$ contains only the prefix embedding of $\rho_j$, implying that $\rho_j$ is a suffix of $\rho_k$ and $\rho_k$ contains exactly these two copies of $\rho_j$.  But $\rho_j$ is also a prefix of $\rho_k$.  Thus $x(\rho_k)=\rho_j$ because if a permutation longer than $\rho_j$ was $x(\rho_k)$, $\rho_{k+1}$ would have more than one copy of $\rho_j$. Since $\rho_k$ contains exactly two copies of $\rho_j$ which are as a prefix and as a suffix, $\rho_j\not\leq i(\rho_k)$. So by Proposition~\ref{outerword}, $C(\rho_{k},\rho_{j})$ is an MSI of $C$.  Thus, by containment minimality, it must be the case that $k=i$.
\hfill\qed

Theorem~\ref{MSIchar} completes the characterization of the MSIs in an interval $[\sigma,\tau]$ of the consecutive pattern poset.  Notice the definitions of the interior and exterior, and the inequality between them,  naturally arise when determining the MSIs.  In particular, they are determined by the lexicographic ordering of the maximal chains.

We are now ready to prove Bernini, Ferrari, and Steingr\'imson's formula using discrete Morse theory.  We have broken the proof up into several cases to make it easier to follow.

\begin{BFST}
In the consecutive pattern poset, if $\sigma\leq \tau$ then 
$$\mu(\sigma,\tau)=\begin{cases}
\mu(\sigma,x(\tau)) & \mbox{if } |\tau|-|\sigma|>2\mbox{ and } \sigma\leq x(\tau)\not\leq i(\tau),\\
1 &\mbox{if } |\tau|-|\sigma|=2\mbox{, } \tau \mbox{ is not monotone, and }\sigma\in\{i(\tau),x(\tau)\},\\
(-1)^{|\tau|-|\sigma|} &\mbox{if } |\tau|-|\sigma|<2,\\
0 &\mbox{otherwise.}\end{cases}$$
\end{BFST}

\prf
Let $[\sigma,\tau]$ be an interval in the consecutive pattern poset.  Suppose first that $|\tau|-|\sigma|<2$. Then $\sigma=\tau$ or $|\sigma|=|\tau|-1$.  By the definition of the M\"obius function, we have $\mu(\sigma,\tau)=1$ in the first case and $\mu(\sigma,\tau)=-1$ in the second case.  Thus, the formula for $\mu(\sigma,\tau)$ holds when $|\tau|-|\sigma|<2$.

Now suppose $|\tau|-|\sigma|=2$.  Then, by the usual M\"obius function recursion, $\mu(\sigma,\tau)=0$ if there is one element in the interval $(\sigma,\tau)$ and $\mu(\sigma,\tau)=1$ when there are 2 elements in the interval $(\sigma,\tau)$.  Since $\tau$ covers at most two elements, these are the only possibilities.  If $\tau$ is monotone, then $\mu(\sigma,\tau)=0$ since $(\sigma,\tau)$ contains a single element. If $\tau$ is not monotone and $\sigma=i(\tau)$, then reducing either the first or last position of $\tau$ gives us an element in $(\sigma,\tau)$, implying $\mu(\sigma,\tau)=1$.  If $\tau$ is not monotone and $\sigma=x(\tau)$, then regardless of whether we reduce the first two positions or last two positions of $\tau$, we get $\sigma$.  Thus, $(\sigma,\tau)$ has 2 elements implying $\mu(\sigma,\tau)=1$.  If the above cases do not hold, then $\sigma$ is either a prefix or a suffix of $\tau$, but not both.  In these cases, $(\sigma,\tau)$ has 1 element implying $\mu(\sigma,\tau)=0$.  Thus, the formula for $\mu(\sigma,\tau)$ holds when $|\tau|-|\sigma|=2$.  

We now turn to the case $|\tau|-|\sigma|>2$.  We will use Theorem~\ref{chainmu} to calculate $\mu(\sigma,\tau)$ from the critical chains in $[\sigma,\tau]$.  By Lemma~\ref{ascent}, the chain id of a critical chain must be decreasing. Since a strong descent is followed by an ascent unless it is the last element in a chain, all the descents must be weak descents except possibly the last one.  It follows that the only maximal chain in $[\sigma,\tau]$ that could be critical is the one which is lexicographically last. Call this chain
$$
C:\tau=\rho_0\stackrel{l_1}{\rightarrow}\rho_1\stackrel{l_2}{\rightarrow}\ldots\stackrel{l_{n-1}}{\rightarrow} \rho_{n-1}\stackrel{l_n}{\rightarrow}\rho_n=\sigma.
$$

Suppose first that $\sigma\leq x(\tau)\not\leq i(\tau)$.  We wish to show that $x(\tau)$ is an element in the chain $C$. Let $k=|\tau|-|x(\tau)|$.  Since $x(\tau)$ is not contained in $i(\tau)$, the prefix of $\tau$ of length $|x(\tau)|+1$ cannot be monotone even if $x(\tau)$ is monotone.  This observation, along with the fact that $\sigma\leq x(\tau)$, allows us to conclude that $|\tau|,|\tau|-1,\ldots,|\tau|-(k-1)$ is a valid beginning for the chain id of a maximal chain $D$ in $[\sigma,\tau]$. Notice that each of these entries is the largest possible entry that does not already appear in the sequence.  Thus, any chain whose chain id differs from the chain id of $D$ in the first $k$ entries is lexicographically earlier than $D$.  So in the chain $C$, $l_1=|\tau|$, $l_2=|\tau|-1$,$\ldots$, $l_k=|\tau|-(k-1)$ which forces $\rho_k=x(\tau)$.

If $\sigma=x(\tau)$, the previous paragraph implies that the sequence $l_1,\ldots,l_k$ is decreasing.  Thus Theorem~\ref{MSIchar} implies $C(\tau,x(\tau))$ is the only interval in $J(C)$.  So by Theorem~\ref{chainmu}, $\mu(\sigma,\tau)=1$.  Of course, in this case $\mu(\sigma,x(\tau))=\mu(\sigma,\sigma)=1$ as well, so the formula holds.

Next we consider $\sigma<x(\tau)$.  Since $l_k$ was the largest possible entry remaining, $l_{k+1}<l_{k}$, implying that $x(\tau)$ is a descent. Let $C'$ be the restriction of $C$ to the interval $[\sigma,x(\tau)]$. We will show that $\#J(C)=2+\#J(C')$, allowing us to apply Theorem~\ref{chainmu} to complete the case $\sigma\leq x(\tau)\not\leq i(t)$.  

Since the sequence $l_1,\ldots,l_{k}$ is decreasing, Theorem~\ref{MSIchar} implies $C(\tau,x(\tau))$ is the first interval in $J(C)$.   We claim that $x(\tau)$ is the second interval in $J(C)$.  If $x(\tau)$ is a strong descent, this follows from the same theorem.  If $x(\tau)$ is a weak descent, $\rho_{k+1}$ is the longest prefix of $x(\tau)$, implying there are at least two copies of $\rho_{k+1}$ contained in $\tau$.  Let $j$ be the the largest value such that $\rho_j$ contains at least two copies of $\rho_{k+1}$.    Since $j$ was picked to be maximum and $\rho_1,\ldots,\rho_{j}$ are weak descents, $\rho_j$ must contain exactly two copies of $\rho_{k+1}$ and they must be the prefix and suffix embeddings.  It follows that $\rho_{k+1}\nleq i(\rho_j)$.  Furthermore, $x(\rho_j)=\rho_{k+1}$ because the prefix with one additional letter, $x(\tau)$, appears only once in $\rho_j$.  Thus, $x(\rho_j)\not\leq i(\rho_j)$.  Since the sequence $l_{j+1},\ldots,l_{k+1}$ is decreasing, Theorem~\ref{MSIchar} implies $C(\rho_j,\rho_{k+1})$ is an MSI in $C$.  By the process of constructing $J(C)$ from $I(C)$, $x(\tau)=C(\rho_j,\rho_{k+1})-C(\rho_0,\rho_k)$ is the second MSI in $J(C)$, proving the claim. Since $x(\tau)$ is an MSI consisting of one element, all the remaining intervals in $J(C)$ are contained in the interval $(\sigma,x(\tau))$.   Therefore, $J(C)=J(C')\cup\{C(\tau,x(\tau)),x(\tau)\}$ and $\#J(C)=2+\#J(C')$.  So by Theorem~\ref{chainmu}, $\mu(\sigma,\tau)=\mu(\sigma,x(\tau))$, proving the formula for $|\tau|-|\sigma|>2$ and $\sigma\leq x(\tau)\not\leq i(\tau)$.

It remains to consider what happens when $|\tau|-|\sigma|>2$ and $\sigma\leq x(\tau)\nleq i(\tau)$ does not hold.  To show $\mu(\sigma,\tau)=0$, we proceed by contradiction.  If $\mu(\sigma,\tau)\neq0$ then, by Theorem~\ref{chainmu}, $J(C)$ must cover $C$. This implies that $J_1=C(\rho_0,\rho_j)$ is an MSI for some $\rho_j$.  Recall that Theorem~\ref{MSIchar} gives two possibilities for MSIs.  If $J_1=\rho_1$ and $\rho_1$ is a strong descent, then since $|\tau|-|\sigma|>2$, $\rho_2$ is an ascent. This contradicts the fact that $C$ has a decreasing chain id.  Alternatively, we must have $\rho_j=x(\tau)$ and, by Theorem~\ref{MSIchar}, $\sigma\leq\rho_j= x(\tau)\nleq i(\tau)$, contradicting our assumption that this  does not hold. So $\mu(\sigma,\tau)=0$, completing the proof. 
\hfill\qed

We can also provide a simple characterization of the homotopy type of an interval in the consecutive pattern poset.  Given two elements $x$ and $y$ of a poset $P$,  the order complex $\Delta(x,y)$ is the abstract simplicial complex whose simplices are the chains in the open interval $(x,y)$.  By using discrete Morse theory to investigate the chains of $(\sigma,\tau)$, one can easily obtain information about the homotopy type of the order complex.

\begin{Thm}
\label{homtype}
Let $[\sigma,\tau]$ be an interval in the consecutive pattern poset.  Then $\Delta(\sigma,\tau)$ is homotopic to a sphere or is contractible, depending on whether the interval has 1 or 0 critical chains, respectively.  In the former case, the dimension of the sphere is $d(C)$, the critical dimension of the critical chain $C$.
\end{Thm}
\prf
In~\cite{BH05}, Babson and Hersh showed that the order complex of a poset is homotopy equivalent to a CW complex with exactly one cell of dimension $d$ for each critical chain of critical dimension $d$ in the interval (as well as a dimension $0$ cell).  In particular, the cell corresponds to a critical simplex in $C$ containing the lowest rank element from each of the intervals in $J(C)$.

In the case of the consecutive pattern poset, Lemma~\ref{ascent} implies that a critical chain must have a decreasing chain id.  Thus, only the lexicographically last chain can be critical.  If $[\sigma,\tau]$ does not have a critical chain, $\Delta(\sigma,\tau)$ is homotopy equivalent to a CW complex with only the $0$-cell.   Thus, the order complex is contractible in this case.    If $[\sigma,\tau]$ has one critical chain $C$, then $\Delta(\sigma,\tau)$ is homotopy equivalent to a CW complex with a $0$-cell and one cell of dimension $d(C)$.  The unique way to attach this cell to the $0$ cell is through a map which is constant on the boundary, resulting in a sphere of dimension $d(C)$.
\hfill\qed

\section{Intervals of Factor Order and the Consecutive Pattern Poset}
\label{cppofo}
The M\"obius function of factor order was determined by Bj\"orner~\cite{aB93}.  We begin this section by defining this partial order.

  Let $A$ be any set.  The \emph{Kleene closure}, $A^*$, is the set of all finite length words over $A$.  So if $w$ is a word and $w(i)$ is the $i^\text{th}$ letter in $w$, then 
  $$A^*=\{w=w(1)w(2)\ldots w(n): 0\le n<\infty \text{ and } w(i)\in A \text{ for all } i\}.$$   The \emph{length} of $w$, denoted $|w|$, is the number of letters in $w$.  \emph{Factor order} on $A^*$ is the partial order on $A^*$ defined by letting $u\leq w$ if $w$ contains a subsequence of consecutive letters $w(i+1)w(i+2)\ldots w(i+n)$ such that $u(j)=w(i+j)$ for $1\leq j\leq n=|u|$.  When $u\leq w$, we call $u$ a \emph{factor} of $w$.  A word $u$ is \emph{flat} if $u(1)=\ldots =u(n)$, where $n=|u|$.  

A \emph{prefix} of a word $w\in A^*$ is a factor of $w$ that includes the first letter of $w$.  Similarly, a \emph{suffix} of $w$ is a factor of $w$ that contains the last letter of $w$.  A prefix or suffix is proper if it is not equal to $w$.  Define the \emph{outer word} $o(w)$ of $w$ to be the longest factor that appears as both a proper prefix and suffix in $w$. Notice that $o(w)$ can be the empty word.  Define the \emph{inner word} $i(w)$ of $w$ to be the factor $i(w)=w(2)...w(n-1)$, where $n=|w|$.

The following theorem of Bj\"orner gives a formula for the M\"obius function of factor order. 

\begin{Thm}[\cite{aB93}]
\label{muofo}
In factor order, if $u\leq w$ then 
$$\mu(u,w)=\begin{cases}
\mu(u,o(w)) & \mbox{if } |w|-|u|>2 \mbox{ and } u\leq o(w)\not\leq i(w),\\
1 &\mbox{if } |w|-|u|=2\mbox{, } w \mbox{ is not flat, and } u\in\{o(w), i(w)\}\\
(-1)^{|w|-|u|} &\mbox{if } |w|-|u|<2,\\
0 &\mbox{otherwise.\hspace*{2.3in}\qed
}\end{cases}$$

\end{Thm}

This formula is remarkably similar to Bernini, Ferrari, and Steingr\'imson's formula.  It turns out that Willenbring's proof of Bj\"orner's formula, which utilizes discrete Morse theory, parallels the one given in the previous section.  It is thus natural to wonder whether there is some relationship between these two posets.  We will give two reasons for this correspondence.  For the first, we will show that certain intervals in $S$ are isomorphic to intervals in factor order over a two letter alphabet.  For the second, we will construct a sequence of posets interpolating between factor order on the positive integers and a poset containing $S$.  All these posets have M\"obius functions satisfying the same recursion and this can be proved uniformly. 

Let $A=\{a,b\}$ and partially order $A^*$ by factor order.  Then it is clear that the number of words of length $n$ is $2^n$.  Consider the set of 213 and 231 avoiding permutations (where these are general, not consecutive, patterns).   Simion and Schmidt~\cite{SS85} were the first to show that the set of such permutations in $S_n$ is counted by $2^{n-1}$.  They gave two proofs, one inductive and one combinatorial.  However, we will use a different combinatorial proof.
Consider constructing a $\sigma=\sigma(1)\sigma(2)\ldots\sigma(n)$ avoiding 213 and 231 by choosing the elements $\sigma(i)$ in increasing order of the index $i$.  Then at each stage, $\sigma(i)$ must be either the largest or the smallest of the elements not already chosen, for if it was not then it would act as the 2 in a copy of either 213 or 231.  So there are two choices for $\sigma(i)$ for all $i<n$ (and only one choice left over for $\sigma(n)$).  Thus, the total number of choices is $2^{n-1}$.

It turns out that the map suggested by the paragraph above,  in which a word in $A^*$ of length $n-1$ is mapped to a 213 and 231 avoiding permutation of length $n$ based on the letter that appears at each position $i$, gives a bijection between these two posets.  Indeed, for the forward map $f$, when we see an $a$ in position $i<n$, we will record the lowest number remaining from the set $\{1,\ldots,n\}$ in position $i$ of the permutation, and when we see a $b$, we will record the highest number remaining from the set $\{1,\ldots,n\}$.  The unused number appears as the last number in the permutation.  Clearly, the inverse map $f^{-1}$ considers each position $i$ of the permutation from left to right (except the last), recording an $a$ if it is the lowest number remaining in the set $\{1,\ldots,n+1\}$, and a $b$ if it is the highest number.

For example, $f(abbab)=165243$ and $f(babab)=615243$, while we have $f^{-1}(1)=\emptyset$, $f^{-1}(12345)=aaaa$, and $f^{-1}(15234)=abaa$.  We now show that this bijection is an order isomorphism.

\begin{Prop}
Let $A=\{a,b\}$.  Then the map $f$ described above is an isomorphism from the poset $A^*$ in factor order to the subposet of the consecutive pattern poset $S$ consisting of all permutations that avoid the patterns $213$ and $231$.
\end{Prop}
\prf
We will only show that $f$ is order preserving as the proof for $f^{-1}$ is similar.  Suppose $u\le w$ in $A^*$ and let $w'$ be the factor of $w$ with $w'=u$.  At each stage of the algorithm for $f$ we always pick the smallest or largest remaining element in constructing $f(w)$.  It follows that $f(w')=f(u)$ is a consecutive pattern in $f(w)$ and so $f(u)\le f(w)$ as desired.
\hfill\qed

\begin{Ques}
For every interval $[\sigma,\tau]$ of the consecutive pattern poset $S$, does there exist an alphabet $A$ and an interval $[u,w]$ in factor order on $A^*$ with $[\sigma,\tau]\cong[u,w]$ (isomorphism of posets)?  What about the same question with the roles of $S$ and $A^*$ reversed?
\end{Ques}

We will now explain the similarity of the two M\"obius recursions by exhibiting infinitely many posets satisfying such a formula.
Let $[0,\infty]$ be the totally ordered set $0<1<2<\cdots<\infty$.  Consider $\bbP^*$ where $\bbP$ is the positive integers.  We are going to construct a partially ordered set $P_k$ on $\bbP^*$ for each $k\in[0,\infty]$ such that $P_0$ will be factor order on $\bbP^*$ and $S$ will be a subposet of $P_\infty$.

We first need to generalize the notion of standardization.  If $w\in\bbP^*$ then define the \emph{height} of $w(i)$ to be
$$
h(w(i))=\#\{j\ |\ \mbox{$1\le j\le w(i)$ and $j\not\in w$}\}.
$$
For example, if $w=24825$ then $h(w(1))=h(w(4))=1$, $h(w(2))=h(w(5))=2$, $h(w(3))=4$.  For $k\in[0,\infty]$ and $w=w(1)w(2)\ldots w(n)$ we define the \emph{$k$-standardization of $w$} to be
$$
\mbox{$\st_k(w)=u(1)u(2)\ldots u(n)$ where $u(i)=w(i)-\min\{h(w(i)),k\}$}
$$
for $1\le i\le n$.  Continuing our example,
\bea
\st_0(24825)&=&24825,\\
\st_1(24825)&=&13714,\\
\st_2(24825)&=&12613,\\
\st_3(24825)&=&12513,\\
\st_k(24825)&=&12413,\\
\eea
for all $k\ge4$.  Note that $\st_0(w)=w$ and $\st_\infty(w)$ is the usual standardization of a sequence of positive integers.  In particular, if the elements of $w$ are distinct then $\st_\infty(w)$ is the standard form of $w$.

To define the poset $P_k$, it suffices to define its covers and then take the transitive closure.  So we let $w=w(1)(2)\ldots w(n)$ cover $\st_k(w(1)\ldots w(n-1))$ and $\st_k(w(2)\ldots w(n))$.  Continuing our example, in $P_2$ the word $w=24825$ would cover $\st_2(2482)=1261$ and $\st_2(4825)=2613$.  This partial order can be given another description as follows.  Call $u$ and \emph{$k$-factor} of $w$ if there is a factor $w'=w(i+1)w(i+2)\ldots w(i+j)$.  Such that 
$$
\st_k^{|w|-|u|}(w')=u.
$$
Then $u\le w$ in $P_k$ if and only if $u$ is a $k$-factor of $w$.  Returning to our example again, a $2$-factor of $w=24825$ is 
$u=\st_2^2(482)=241$ and in $P_2$ we have the interval $[u,w]=\{241,1261,2613,24825\}$.  

It should be clear from the definitions that $P_0$ is factor order on $\bbP^*$, and that the consecutive pattern poset $S$ is the subposet of $P_\infty$ consisting of all words which are permutations.  It is also easy to see that $S$ is  a \emph{convex} subposet of $P_\infty$  in that if $u,w\in S$ then all elements of the interval $[u,w]$ in $P_\infty$ are in $S$.  It follows that $\mu(u,w)$ is the same in $P_\infty$ and in $S$.

To state the recursive formula for the M\"obius function of $P_k$, we need the analogues of interior/inner and exterior/outer.  These should come as no supprise.  Define the \emph{$k$-interior} of $w=w(1)\ldots w(n)$ to be
$$
i_k(w)=\st_k^2(w(2)\ldots w(n-1)).
$$
A \emph{$k$-prefix} of $w$ is a $k$-factor where the subword standardized contains $w(1)$ and similarly for a \emph{$k$-suffix}.  The \emph{$k$-exterior} of $w$ is $x_k(w)$, the longest proper $k$-factor which is both a $k$-prefix and a $k$-suffix.

The last ingredient of the recursive formula for $\mu$ in $P_k$ is the analogue of a monotone permutation in $S$ and of a flat word in factor order.  Say that $w\in P_k$ is \emph{$k$-irreducible} if it covers only one element in $P_k$.  The reason for this terminology is that, in a poset which is a latteice, an element is join irreducible if it covers only one element.  To see why this is the correct notion to capture monotonicity and flatness, we will prove a lemma.  In it, a word $w$ of length $n$ is \emph{monotone increasing} if $w(1)<w(2)<\cdots<w(n)$ and similarly for monotone decreasing.
\begin{Lem}
\label{k-irr}
In $P_k$ a word $w=w(1) w(2)\ldots w(n$), $n\ge1$, is $k$-irreducible if and only 
\begin{enumerate}
\item $w$ is monotone increasing and $w(i)\le k+i-1$ for $1\le i\le n$, or
\item $w(1)=w(2)=\cdots=w(n)$, or
\item $w$ is monotone decreasing and $w(i)\le n-i+k$ for $1\le i\le n$.
\end{enumerate}
\end{Lem}
\prf
It is a simple matter to check that any word satisfying one of the three given conditions covers a unique element, so we will only prove the forward direction.  Let 
$$
u=\st_k(w(1)\ldots w(n-1),\ v=\st_k(w(2)\ldots w(n))
$$
Assume that $w$ is $k$-irreducible so that $u=v$.

Suppose first that $w(1)<w(2)$.  It follows that $h(w(1))\le h(w(2))$ (where height is being computed to obtain $u$).  This implies $u(1)<u(2)$.  But $u=v$, so $v(1)<v(2)$ which gives $w(2)<w(3)$.  Continuing in this way, we obtain that $w$ is monotone increasing.  Assuming $w(1)>w(2)$ or $w(1)=w(2)$ lead to the monotone decreasing or all equal cases, respectively.

From here on, we will assume that $w$ is monotone increasing as the third case is similar and the second is easier.  Also assume, towards a contradiction, that $w(1)>k$.  Since $w$ is increasing, $h(w(1))=w(1)-1\ge k$.  So $u(1)=w(1)-k$.  Applying the same reasoning to $w(2)>k+1$ when producing $v$ gives $v(1)=w(2)-k$.  But then $u(1)=w(1)-k<w(2)-k=v(1)$, contradicting the fact that $u=v$.  If follows that $w(1)\le k$.

We now show that $w(1)\le k$ implies $w(2)\le k+1$ as the same argument can be used to inductively show that $w(i)\le k+i-1$ for all $1\le i\le n$.  Indeed, $w(1)\le k$ gives $u(1)=1$ and so $v(1)=1$ as well. But this forces $w(2)\le k+1$ from the computation for $v$.  This finishes the proof.
\hfill\qed

Note that in $P_0$, the first case can not happen as it would require $w(1)\le0+1-1=0$.  Similarly, the third case can not happen by considering $w(n)$.  So $0$-irreducible corresponds to flat.  In $S$, case 2 can not happen by definition.  Furthermore, since $k=\infty$, the restrictions on $w(i)$ hold trivially.  So we are left with only the monotone condition.

\begin{Thm}
In $P_k$, $k\in[0,\infty]$, if $u\le w$ then
$$\mu(u,w)=\begin{cases}
\mu(u,x_k(w)) & \mbox{if } |w|-|u|>2\mbox{ and } u\leq x_k(w)\not\leq i_k(w),\\
1 &\mbox{if } |w|-|u|=2\mbox{, } w \mbox{ is not $k$-irreducible, and }u\in\{i_k(w),x_k(w)\},\\
(-1)^{|w|-|u|} &\mbox{if } |w|-|u|<2,\\
0 &\mbox{otherwise.}\end{cases}$$
Furthernore, $\Delta(u,w)$  is homotopic to a sphere or is contractible, depending on whether the interval has 1 or 0 critical chains, respectively.  In the former case, the dimension of the sphere is $d(C)$, the critical dimesion of the critical chain $C$.\hfill\qed
\end{Thm}

The proof of this theorem is virtually identical to the proofs of Theorems~\ref{mucpp} and~\ref{homtype} given previously.  One need only use either the definition of $k$-irreducible or Lemma~\ref{k-irr} as appropriate every time a fact about monotone permutations is needed.  So it turns out that consecutive pattern order and factor order are just two extremes of an infinite family of posets whose M\"obius functions have the same structure.


\begin{thebibliography}{1}

\bibitem{BH05}
{\sc Babson, E., and Hersh, P.}
\newblock Discrete {M}orse functions from lexicographic orders.
\newblock {\em Trans. Amer. Math. Soc. 357}, 2 (2005), 509--534 (electronic).

\bibitem{BFS11}
{\sc Bernini, A., Ferrari, L., and Steingr\'imson, E.}
\newblock The {M}\"obius function of the consecutive pattern poset.
\newblock Preprint {\texttt{arXiv:1103.0173}}.

\bibitem{aB93}
{\sc Bj{\"o}rner, A.}
\newblock The {M}\"obius function of factor order.
\newblock {\em Theoret. Comput. Sci. 117}, 1-2 (1993), 91--98.
\newblock Conference on Formal Power Series and Algebraic Combinatorics
  (Bordeaux, 1991).

\bibitem{rF95}
{\sc Forman, R.}
\newblock A discrete {M}orse theory for cell complexes.
\newblock In {\em Geometry, topology, \& physics}, Conf. Proc. Lecture Notes
  Geom. Topology, IV. Int. Press, Cambridge, MA, 1995, pp.~112--125.

\bibitem{rF02}
{\sc Forman, R.}
\newblock A user's guide to discrete {M}orse theory.
\newblock {\em S\'em. Lothar. Combin. 48\/} (2002), Art.\ B48c, 35 pp.
  (electronic).

\bibitem{SS85}
{\sc Simion, R., and Schmidt, F.~W.}
\newblock Restricted permutations.
\newblock {\em European J. Combin 6}, 4 (1985), 383--406.

\bibitem{RW11}
{\sc Willenbring, R.}
\newblock The {M}\"obius function of generalized factor order.
\newblock Submitted.

\end{thebibliography}
\end{document}